\documentclass{amsart}
\usepackage[centertags]{amsmath}
\usepackage{amsfonts}
\usepackage{amssymb}
\usepackage{amsthm}
\usepackage{graphicx}
\usepackage{caption}
\usepackage{subcaption}
\usepackage[all]{xy}
\usepackage{tikz}
\usepackage{eepic}
\usepackage{curve2e}

\newtheorem{thm}{Theorem}
\newtheorem{lem}[thm]{Lemma}
\newtheorem{conj}[thm]{Conjecture}

\def\dl#1{{\displaystyle{#1} }}

\def\ints{{{\mathbb Z}}}
\def\Ex{{ \hbox{Ext}_A^{s,t} (} }
\def\Ext#1{{ \hbox{Ext}_{#1}} }

\def\z2{{(\ints /2, \ints /2) }}
\def\f2{{(\F_2,\F_2)}}
\def\la{{\langle}}
\def\ra{{\rangle}}
\def\ASS{{Adams spectral sequence}}
\def\SS{{spectral sequence}}

\def\mapright#1{{\smash{\mathop{\longrightarrow}\limits^{#1}}}}
\def\mapleft#1{{\smash{\mathop{\longleftarrow}\limits^{#1}}}}\def\elt{\circle*{3}}
\def\F{{\mathbb F}}
\newcommand{\sslash}{\mathbin{/\mkern-5mu/}}

\def\elt{\circle*{3}}
\def\hoa{\line(0,1){20}}

\def\dia{\vector(-1,1){18}}

\def\tower{\vector(0,1){20}}

\def\scsz{\scriptsize}

\usetikzlibrary{matrix}

\global\overfullrule 0pt 
\begin{document}

\title[Mahowald's Conjecture on $\Ext{A(n)}(\F_2,\F_2)$]{On a Conjecture of Mahowald on the Cohomology of Finite
Sub-Hopf algebras of the Steenrod Algebra}
\author{Paul Shick}
\address{John Carroll University, University Hts OH 44118}
\maketitle

\centerline{\today}

\bigskip

  Mahowald's conjecture arose as part of a program attempting to view chromatic phenomena in stable homotopy theory through the lens of the classical Adams spectral sequence. The conjecture predicts the existence of nonzero classes in the cohomology of the finite sub-Hopf algebras $A(n)$ of the mod 2 Steenrod algebra that correspond to generators in the homotopy rings of certain periodic spectra.  The purpose of this note is to present a proof of the conjecture.

\smallskip
\section{Introduction and Statement of Results} \label{Intro}

  To provide some context for the conjecture, here's a
brief summary of what has been known about how elements that detect periodic phenomena
appear in the cohomology of finite Hopf-subalgebras of the Steenrod
algebra.  For a Hopf-subalgebra $B$ of the mod $p$ Steenrod algebra, we'll often use the notation $\dl{H^*(B)}$ as an abbreviation for the cohomology of $B,$ $\dl{\Ext{B}(\F_p,\F_p).}$

\smallskip
Let $A(n)$ denote the finite subHopf algebra of the Steenrod
algebra generated by
$$Sq^0,Sq^1,\dots,Sq^{2^n}$$
if $p=2$ and by
$$\beta,P^1,\dots,P^{p^{n-1}},$$
if $p$ odd.  Let $\dl{E(n) = E(Q_0,Q_1,\dots,Q_n)}$ denote the
subalgebra of $A(n)$ generated by the Milnor generators.  Since
the Johnson-Wilson spectrum $BP\la n \ra$ has  $\dl{H^*(BP\langle
n \rangle; \F_p) = A\sslash E(n)}$, we can calculate $\dl{\pi_*(BP \langle n
\rangle)}$ from the classical {\ASS}:
\begin{eqnarray*}
 E_2^{s,t} & = & \Ex H^*(BP \langle n \rangle),\F_p)\cr
& = & \Ex A\sslash E(n),\F_p) \\
& \cong  & \hbox{Ext}_{E(n)}(\F_p,\F_p)\\
& \cong  & \F_p[q_0,q_1,\dots,q_n].
\end{eqnarray*}
We observe that $E_2 = E_\infty$, since the generators are concentrated in
even degrees, and
$$ \pi_*(BP \langle n \rangle) \cong
\ints_{(p)}[v_1,\dots ,v_n],$$ where $|v_i| = 2p^i -2$ (and the generators are in Adams
filtration 1).  We hereafter denote the generators of
$\dl{\Ext{E(n)}(\F_p,\F_p)}$ by $v_i$, for $i \ge 0.$
 The inclusion  \hfill \break
 $i:E(n) \hookrightarrow A(n)$ induces the restriction
homomorphism in cohomology
$$i^* : \hbox{Ext}_{A(n)}(\F_p,\F_p) \to
\hbox{Ext}_{E(n)}(\F_p,\F_p).$$

\smallskip
When one attempts to understand chromatic phenomena via the classical Adams spectral sequence,
odd primes are easier to handle than the $p=2$ case, as demonstrated by the following result from \cite{S} (although work on the prime 2 case was done first).
\smallskip

\begin{thm}  For $p$ any odd prime, there are classes defined and nonzero in \break \hfill $\dl{\hbox{Ext}_{A(n)}(\F_p,\F_p)}$
that form a polynomial
subalgebra:
\[
\F_p[v_0,v_1^{p^n},v_2^{p^{n-1}},\dots,v_n^p] \subset
\hbox{Ext}_{A(n)}(\F_p,\F_p),
\]
where the generators restrict to the obvious classes in
$H^*(E(n)).$
\end{thm}

The proof uses a careful analysis of the Cartan-Eilenberg {\SS}
for the extension of Hopf algebras
$$ {\Bbb F}_p \to P(n) \to A(n)_* \to D(n) \to {\Bbb F}_p$$
where $P(n)$ is the truncated polynomial algebra on the $\xi_i$s
and $D(n)$ is the exterior algebra on the $\tau_i$s.
 The spectral sequence collapses for
odd primes, as one can see by filtering the dual Steenrod algebra by the number of $\tau$s.
This allows one to move directly from understanding the coaction on the cohomology of $D(n)$
to seeing the appropriate classes in $\dl{\Ext{A(n)}(\F_p,\F_p).}$

\bigskip
For the prime 2, the situation is more difficult, but there were partial results available (\cite{MS1}):

\begin{thm} There are classes defined and nonzero in $\dl{\hbox{Ext}_{A(n)}\f2}$ that form a polynomial
subalgebra
\[
\F_2[v_0,v_1^{N_1},v_2^{N_2},\dots,v_n^{N_n}] \subset
\hbox{Ext}_{A(n)}\f2,
\]
 where the generators restrict to the
obvious classes in $H^*(E(n)).$
\end{thm}

The proof uses results of Lin \cite{Lin} and Wilkerson
\cite{Wilk} that show that the restriction homomorphism $i^*$ is
onto in infinitely many positive degrees. Note that this argument defines these
generators $\dl{v_i^{N_i}}$ only as cosets for $p=2$ and
$1 \le i \le n-1.$  The top class was explicitly identified as $\dl{v_n^{2^{n+1}}}$  and shown to be a non-zero divisor in the cohomology ring $\dl{\Ext{A(n)}\f2 }$, using a {\SS}
based on a Koszul-type resolution of $A(n)\sslash A(n-1)$. This spectral sequence first appeared in \cite{DM1} and was referred to as the ``Koszul spectral sequence" in \cite{MS1}. It has recently been rechristened as the Davis-Mahowald spectral sequence by Bruner, Rognes and their coworkers in \cite{BR} and \cite{Rog}.
\smallskip

The exponents of the lower $\dl{v_i}$s have been more mysterious.  Low dimensional calculations led to the following
conjecture, the proof of which is the goal of this note.

\begin{conj} \label{Conj}  For any natural numbers $n$ and $k,$ the class $\dl{ v_n^{2^{n+1+k}}}$ is defined and nonzero in $\dl{\Ext{A(j)}\f2}$ for all $j = n, n+1, \dots, 2n+k.$ \end{conj}

\bigskip

This conjecture is originally due to Mahowald around 1980,
although it first appeared in print in \cite{MS1}.  Note that the conjecture predicts that $\dl{v_n^{2^{n+1+k}}}$ is nonzero in the cohomology of $A(j)$ for the largest possible $j$:
If $\dl{v_n^{2^{n+1+k}}}$ is defined and nonzero in $\dl{\Ext{A(j)}\f2}$ for $j > 2n+k,$ then the class would persist to $\dl{\Ext{A}\f2}$ by the Adams Approximation Theorem \cite{Adams}. This class would then be a permanent cycle in the classical Adams spectral sequence for the homotopy ring of the (2-complete) sphere, persisting to a nonnilpotent element in $\dl{\pi_0((S^0)^\wedge)}$, contradicting Nishida's Theorem \cite{Nishida}.

\smallskip
For those who prefer lists, the main theorem can be restated as:
For any natural number $n,$ we have

\[\F_2[v_0,v_1^{2^{2n}},v_2^{2^{2n-1}},v_3^{2^{2n-2}},\dots,v_n^{2^{n+1}}, v_{n+1}^{2^{n+2}},
v_{n+2}^{2^{n+3}}, \dots , v_{2n}^{2^{2n+1}}] \subset
\hbox{Ext}_{A(2n)}\f2
\]

and
\[\F_2[v_0,v_1^{2^{2n+1}},v_2^{2^{2n}},v_3^{2^{2n-1}},\dots,v_n^{2^{n+2}}, v_{n+1}^{2^{n+2}},
v_{n+2}^{2^{n+3}}, \dots , v_{2n+1}^{2^{2n+2}}] \subset
\hbox{Ext}_{A(2n+1)}\f2 .
\]

\smallskip

 The proof of the conjecture
presented here is inductive, using as its main tools the Davis-Mahowald spectral sequence and the classical May spectral sequence (augmented by Nakamura's squaring operations \cite{Nak}).
\smallskip

A very different approach to the conjecture can be found in chapter 5 of Singer's monograph \cite{Sing}.  He considers the Cartan-Eilenberg {\SS} for the 2-primary extension of Hopf algebras
\[ E(n) \to A(n) \to {\mathcal D}A(n-1),\]
where $\dl{{\mathcal D}A(n-1)}$ has algebra and coalgebra structures identical to $A(n-1)$ but with the gradings doubled.  This Cartan-Eilenberg spectral sequence for the cohomology of $A(n)$ fails to collapse for the prime 2, so one needs to look closely at the (right) action of $\dl{{\mathcal D}A(n-1)}$ on $\dl{\hbox{Ext}_{E(n)}\f2 \cong \F_2[v_0,v_1,v_2,\dots]: }$

\[
    v_k  D({\overline{Sq}}^{2^i}) =
    \begin{cases}
        v_{k-1} & \text{if } k=i-1;\\
        0 & \text{otherwise. }
           \end{cases}
  \]
 One then needs to compute the image of the restriction homomorphism
 \[
\Ext{A(n)}(\F_2,\F_2) \to \left( \Ext{E(n)}(\F_2,\F_2) \right)^{{\mathcal D}\overline{A}(n-1)} .
\]
Singer observes that the lowest power of $\dl{v_k}$ invariant under the coaction of $\dl{{\mathcal D}\overline{A}(n-1)}$ is
$\dl{v_k^{2^{n+1-k}},}$ exactly the lowest power predicted to live in the cohomology of $A(n)$ by Conjecture \ref{Conj}.  He suggests that a careful analysis of this Cartan-Eilenberg spectral sequence will likely shed light on the conjecture and other questions about the cohomology of finite Hopf-subalgebras of the mod 2 Steenrod algebra.
\bigskip\bigskip

\noindent It is a pleasure to acknowledge helpful suggestions from John Rognes and Doug Ravenel that have improved the exposition. The author also thanks the referee for his/her careful reading and detailed recommendations. He is especially grateful to Mark Mahowald, whose profound insights opened up this area of inquiry.

\bigskip
\bigskip

\section{The Davis-Mahowald Spectral Sequence} \label{DMSS}

\smallskip
The two main tools used in the proof of the conjecture are the now-classical May spectral sequence and the spectral sequence developed by Davis and Mahowald in their study of $\dl{\Ext{A(2)}(P_k,\F_2),}$ \cite{DM1}.  Because the May spectral sequence is so ubiquitous, there's no need to include details about its construction here.  We note, however, that we'll use the more modern notation for classes in the May {\SS}, found in such sources as \cite{Ravbook}.  We provide some basic information on the construction of the Davis-Mahowald spectral sequence, based on their
variant of the Koszul resolution.  The following material is
based on presentations in \cite{DM1} and \cite{MS1}, augmented by Rognes' later work in \cite{Rog}.
\smallskip

The idea, originally due to Davis and Mahowald, is to use a sort of
``sideways" version of the traditional Koszul resolution to allow
one to compute $\Ext{A(n)}(M,\F_2)$ if one understands
$\Ext{A(n-1)}(M \otimes N ,\F_2)$ for certain $A(n-1)$-modules $N.$  We begin by observing that (at the prime 2) the dual Steenrod algebra
\[
A_* \cong \F_2 [\xi_1, \xi_2 , \dots ]
\]
is a module over $A$, with the action given by the total squaring
operation
\[ Sq \xi_i = \xi_i + \xi_{i-1}^2.
\]
Observe that
\[
\left( A(n) \otimes_{A(n-1)}\F_2 \right)^* \cong  E(\zeta_1^{2^n}, \zeta_2^{2^{n-1}},\dots,\zeta_{n+1}),
\]
 as right $A(n)$-modules, where the action on the exterior algebra is given by
\[
Sq^{2^k} \zeta_{n-j+1}^{2^k}  = \zeta_{n-j}^{2^{k+1}}
\]
and
\[
Sq^{2^n} \zeta_{1}^{2^n}  = 1,
\]
extended by the Cartan formula, where $\zeta_i = \chi(\xi_i)$.  We note that $\dl{E(\zeta_1^{2^n}, \zeta_2^{2^{n-1}},\dots,\zeta_{n+1})}$ is an
$A(n)$-module, but not an $A(n+1)$-module, because of the Adem
relations.
For each $n$, let
\[
R=R(n) =\F_2[\zeta_1^{2^n}, \zeta_2^{2^{n-1}},\dots,\zeta_{n+1}]
\]
with the same $A$-action on the generators, extended by the Cartan formula.  For each $\sigma \ge 0,$ let
$\dl{R^\sigma=R^\sigma(n)}$ denote the sub vector space of $R$
spanned by homogeneous polynomials of degree $\sigma.$  Here are pictures of the first few $\dl{R^i(n})$s:

\begin{figure}[h]
\[
\xymatrix@R-2em@C-2em{
\zeta_2 & \circ \ar@{-}[d]^{Sq^1}
\\
\zeta_1^2 & \circ  \\
}
\]
\caption{$\dl{R^1(1)}$}\label{fig:R11}
\end{figure}
\begin{figure}[h]
\[
\xymatrix@R-2em@C-2em{
\zeta_3 & \circ \ar@{-}[d]^{Sq^1}
\\
\zeta_2^2 & \circ \ar@/^1pc/@{-}[dd]^{Sq^2}
\\ \\
\zeta_1^4 & \circ
}
\]
\caption{$\dl{R^1(2)}$}\label{fig:R12}
\end{figure}
\begin{figure}[h]
\[
\xymatrix@R-2em@C-2em{
\zeta_3^2 & \circ \ar@/^1pc/@{-}[dd]^{Sq^2}
\\
\zeta_2^2 \zeta_3 & \circ \ar@{-}[d]^{Sq^1} \ar@/_1pc/@{-}[dd]
\\
\zeta_2^4 & \circ
\\
\zeta_1^4 \zeta_3 & \circ \ar@{-}[d]^{Sq^1}
\\
\zeta_1^4 \zeta_2^2 & \circ \ar@/^1pc/@{-}[dd]^{Sq^2}
\\ \\
\zeta_1^8 & \circ
}
\]
\caption{$\dl{R^2(2)}$}\label{fig:R22}
\end{figure}
\begin{figure}[h]
\[
\xymatrix@R-2em@C-2em{
\zeta_4 & \circ \ar@{-}[d]^{Sq^1}
\\
\zeta_3^2 & \circ \ar@/^1pc/@{-}[dd]^{Sq^2}
\\ \\
\zeta_2^4 & \circ \ar@{-} `r[dddd] `[dddd]^{Sq^4} [dddd]
\\ \\ \\ \\
\zeta_1^8 & \circ
}
\]
\caption{$\dl{R^1(3)}$}\label{fig:R13}
\end{figure}

\break
For future reference, we also include diagrams of the exterior algebras:
\begin{figure}[h]
\[
\xymatrix@R-2em@C-2em{
\zeta_1^2\zeta_2 & \circ \ar@/^1pc/@{-}[dd]^{Sq^2}
\\ \\
\zeta_2 & \circ \ar@{-}[d]^{Sq^1}
\\
\zeta_1^2 & \circ \ar@/^1pc/@{-}[dd]^{Sq^2}
\\ \\
1 & \circ
}
\]
\caption{$\dl{(A(1)\sslash A(0))_*}$}\label{fig:R14}
\end{figure}

\begin{figure}[h]
\[
\xymatrix@R-2em@C-2em{
\zeta_1^4 \zeta_2^2 \zeta_2 & \circ \ar@{-} `r[dddd] `[dddd]^{Sq^4} [dddd]
\\ \\ \\ \\
\zeta_2^2 \zeta_3 & \circ \ar@/^1pc/@{-}[dd]^{Sq^2}
\\ \\
\zeta_1^4 \zeta_3 & \circ \ar@{-}[d]^{Sq^1} \ar@{-} `l[dddd] `[dddd]^{Sq^4} [dddd]
\\
\zeta_1^4 \zeta_2^2 & \circ \ar@{-} `r[dddd] `[dddd]^{Sq^4} [dddd]
\\ \\ \\
\zeta_3 & \circ \ar@{-}[d]^{Sq^1}
\\
\zeta_2^2 & \circ \ar@/^1pc/@{-}[dd]^{Sq^2}
\\ \\
\zeta_1^4 & \circ \ar@{-} `r[dddd] `[dddd]^{Sq^4} [dddd]
\\ \\ \\ \\
1& \circ
}
\]
\caption{$\dl{(A(2)\sslash A(1))_*}$}\label{fig:R14}
\end{figure}

\break

Then for every $n$ we
have the following sequence:
\[
0 \to \F_2 \to \left(A(n) \otimes_{A(n-1)}\F_2 \right)^* \otimes R^0(n) \to \left(A(n) \otimes_{A(n-1)}\F_2 \right)^* \otimes R^1(n) \to \cdots
\]
with the maps given by
\[
\delta[(x_1x_2\dots x_r) \otimes p] = \sum_{j=1}^r (x_1 \dots
\widehat{x_j}
    \dots x_r) \otimes x_jp.
\]
This sequence is exact by the standard Koszul resolution
argument (See \cite{HS} p 243 or \cite{DM1}.) and is dual to
\begin{eqnarray}\label{KR}
\,\,\,\,\,\,  0 \longleftarrow \F_2 \longleftarrow A(n) \otimes_{A(n-1)} (R^0(n))^*
    \longleftarrow A(n) \otimes_{A(n-1)} (R^1(n))^* \longleftarrow \cdots
\end{eqnarray}
Given an $A(n)$-module $M$, apply the functor
$\dl{\Ext{A(n)}(\,-\,,\F_2)}$ to the resolution above, tensored with
$M$.  We obtain a \SS
\[
E_1^{\sigma,s,t} = \hbox{Ext}_{A(n-1)}^{s-\sigma,t}((R^\sigma(n))^*
\otimes M, \F_2)
    \Longrightarrow \hbox{Ext}_{A(n)}^{s,t}(M,\F_2).
\]
Note that the homomorphisms in the resolution are given by
$A(n)$-module maps which are {\underline{not}} extended $A(n-1)$-module maps.
Thus $d_1$-differentials in the Davis-Mahowald {\SS} are induced by the
$\dl{Sq^{2^n}}$-action in $R^\sigma(n)^* \otimes M$.

\smallskip As an (easy and familiar) example, we'll use this Davis-Mahowald {\SS} to compute
$\dl{\hbox{Ext}_{A(1)}\f2},$ the $E_2$ term of the classical {\ASS}
for $bo_*$.  Here we use as our starting point the fact that the
dual of  $A(1)\sslash A(0)$ is isomorphic to $\dl{ E(\xi_1^2,\xi_2)}$
as $A(1)$-coalgebras.  Note that the gradings are quite important
here: $\xi_1^2$ is in dimension 2 and $\xi_2$ is in dimension 3.
We set up the Koszul resolution of (\ref{KR}) in the $n=1$ case,
then apply the functor $\dl{\hbox{Ext}_{A(1)}^{*,*}(\,-\,,\F_2)}$,
obtaining the Davis-Mahowald {\SS}.  Here the $E_1$ term consists of
$\dl{\hbox{Ext}_{A(0)}} $s, with the $\dl{d_1}$-differentials induced by the
$Sq^2$-action on the $\dl{R^\sigma}$s (not the $\dl{Sq^1}$ action).  A
picture of the $\dl{E_1}$ term follows, in the standard $(t-s,s)$
chart form, with the $\sigma$-filtration of the classes labeled
appropriately.

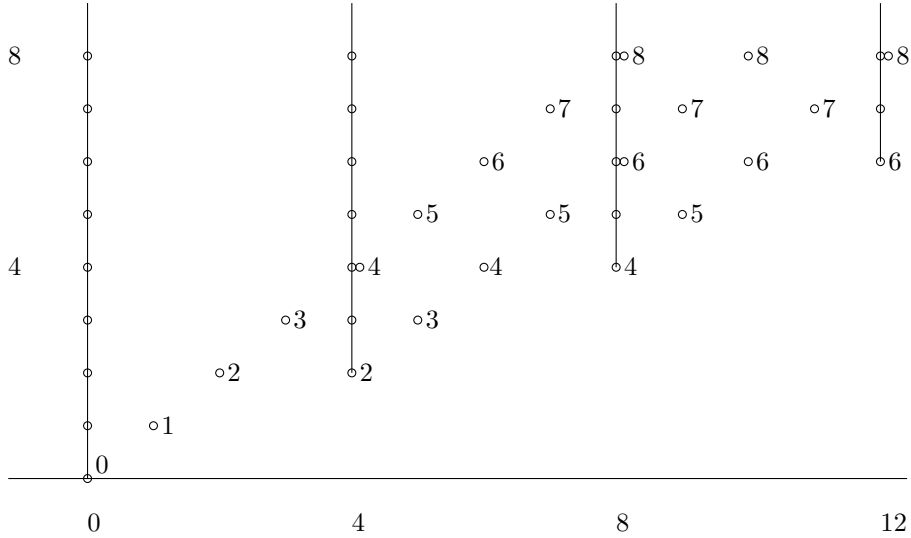
\begin{center}
\begin{picture}(360,200) \label{KSSA1}
\put(0,20){\line(1,0){340}} \put(130,0){$4$} \put(230,0){$8$}
\put(330,0){$12$} \put(0,97){$4$} \put(0,177){$8$} \put(30,0){$0$}
\put(30,20){\line(0,1){180}} \put(30,20){\elt} \put(30,40){\elt}
\put(30,60){\elt} \put(30,80){\elt} \put(30,100){\elt}
\put(30,120){\elt} \put(30,140){\elt} \put(30,160){\elt}
\put(30,180){\elt} \put(33,22){$0$} \put(55,40){\elt}
\put(58,37){$1$} \put(80,60){\elt} \put(83,57){$2$}
\put(130,60){\elt} \put(130,60){\line(0,1){140}} \put(133,57){$2$}
\put(130,80){\elt} \put(130,100){\elt} \put(130,120){\elt}
\put(130,140){\elt} \put(130,160){\elt} \put(130,180){\elt}
\put(105,80){\elt} \put(108,77){$3$} \put(155,80){\elt}
\put(158,77){$3$} \put(133,100){\elt} \put(136,97){$4$}
\put(180,100){\elt} \put(182,97){$4$}
\put(230,100){\line(0,1){100}} \put(230,100){\elt}
\put(230,120){\elt} \put(230,140){\elt} \put(230,160){\elt}
\put(230,180){\elt} \put(233,97){$4$} \put(155,120){\elt}
\put(158,117){$5$} \put(205,120){\elt} \put(208,117){$5$}
\put(255,120){\elt} \put(258,117){$5$} \put(180,140){\elt}
\put(183,137){$6$} \put(233,140){\elt} \put(236,137){$6$}
\put(280,140){\elt} \put(283,137){$6$}
\put(330,140){\line(0,1){60}} \put(330,140){\elt}
\put(330,160){\elt} \put(330,180){\elt} \put(333,137){$6$}
\put(205,160){\elt} \put(208,157){$7$} \put(255,160){\elt}
\put(258,157){$7$} \put(305,160){\elt} \put(308,157){$7$}
\put(233,180){\elt} \put(236,177){$8$} \put(280,180){\elt}
\put(283,177){$8$} \put(333,180){\elt} \put(336,177){$8$}
\end{picture}
\end{center}

Now we need to sort out the differentials in the Davis-Mahowald \SS.  Since
$\dl{d_1: E_1^{\sigma,s,t} \to E_1^{\sigma +1,s+2,t}}$ (i.e., one up
and one to the left on the chart), we see a potential $\dl{d_1}$ from
the second filtration 2 class to the first filtration 3 class. To
check whether or not this is nonzero, we need to look at the
homomorphism
\[
A(1)\sslash A(0) \otimes R^2(1)^* \mapleft{\delta_2^*} A(1)\sslash A(0) \otimes
R^3(1)^*,
\]
dual to
\[
E(1) \otimes R^2 \mapright{\delta_2} E(1) \otimes R^3.
\]
For any class $\dl{p \in R^2}$, we see that
\begin{eqnarray*}
\delta_2(1\otimes p) = & 0 \\
\delta_2(\zeta_1^2\otimes p)  = & 1 \otimes \zeta_1^2p \\
\delta_2(\zeta_2\otimes p)  = & 1 \otimes \zeta_2p \\
\delta_2(\zeta_1^2\zeta_2\otimes p)  = & \zeta_1^2 \otimes
\zeta_2p
        + \zeta_2 \otimes \zeta_1^2 p. \\
        \end{eqnarray*}
To see which polynomials $p$ yield a nonzero $\dl{d_1}$, we look at
the action of $\dl{Sq^2}$ on $p$.  Since $\dl{\zeta_1^4 Sq^2 =0}$ and
$\zeta_1^2 \zeta_2 Sq^2 =0$, the classes spawned by these in Ext
will have zero $d_1$s.  But $\zeta_2^2 Sq^2 = \zeta_1^4,$ so that
the class this gives in Ext (the second filtration 2 class) must
hit the first filtration 3 class.
\smallskip

This sort of reasoning establishes all of the $\dl{d_1}$s in the
following picture:

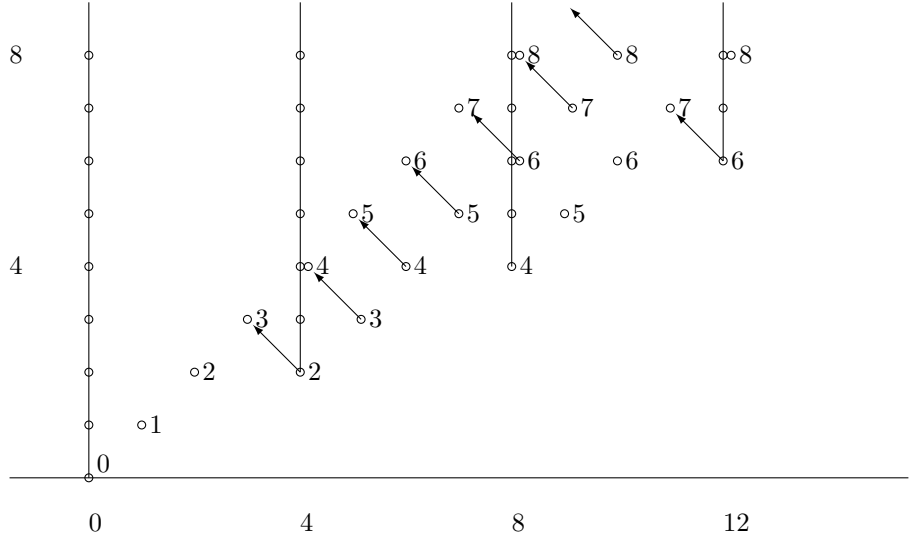
\begin{center}
\begin{picture}(300,200) \label{KSSA1dif}
\put(0,20){\line(1,0){340}} \put(110,0){$4$} \put(190,0){$8$}
\put(270,0){$12$} \put(0,97){$4$} \put(0,177){$8$} \put(30,0){$0$}
\put(30,20){\line(0,1){180}} \put(30,20){\elt} \put(30,40){\elt}
\put(30,60){\elt} \put(30,80){\elt} \put(30,100){\elt}
\put(30,120){\elt} \put(30,140){\elt} \put(30,160){\elt}
\put(30,180){\elt} \put(33,22){$0$} \put(50,40){\elt}
\put(53,37){$1$} \put(70,60){\elt} \put(73,57){$2$}
\put(110,60){\elt} \put(110,60){\line(0,1){140}} \put(113,57){$2$}
\put(110,60){\dia} \put(110,80){\elt}
\put(110,100){\elt} \put(110,120){\elt} \put(110,140){\elt}
\put(110,160){\elt} \put(110,180){\elt} \put(90,80){\elt}
\put(93,77){$3$} \put(133,80){\elt} \put(136,77){$3$}
\put(133,80){\dia} \put(113,100){\elt}
\put(116,97){$4$} \put(150,100){\elt} \put(153,97){$4$}
\put(150,100){\dia} \put(190,100){\line(0,1){100}}
\put(190,100){\elt} \put(190,120){\elt} \put(190,140){\elt}
\put(190,160){\elt} \put(190,180){\elt} \put(193,97){$4$}
\put(130,120){\elt} \put(133,117){$5$} \put(170,120){\elt}
\put(173,117){$5$} \put(170,120){\dia}
\put(210,120){\elt} \put(213,117){$5$} \put(150,140){\elt}
\put(153,137){$6$} \put(193,140){\elt} \put(196,137){$6$}
\put(193,140){\dia} \put(230,140){\elt}
\put(233,137){$6$} \put(270,140){\line(0,1){60}}
\put(270,140){\elt} \put(270,160){\elt} \put(270,180){\elt}
\put(273,137){$6$} \put(270,140){\dia}
\put(170,160){\elt} \put(173,157){$7$} \put(213,160){\elt}
\put(216,157){$7$} \put(213,160){\dia}
\put(250,160){\elt} \put(253,157){$7$} \put(193,180){\elt}
\put(196,177){$8$} \put(230,180){\elt} \put(233,177){$8$}
\put(230,180){\dia} \put(273,180){\elt}
\put(276,177){$8$}
\end{picture}
\end{center}

This leaves the expected picture for $\dl{\hbox{Ext}_{A(1)}\f2}$.

\bigskip

The Davis-Mahowald {\SS} gives an easy proof that $\dl{v_n^{2^{n+1}}}$ is defined and nonzero in $\dl{\hbox{Ext}_{A(n)}^{2^{n+1},2^{n+1}(2^{n+1}-1)}\f2},$ which we outline here.  We note that the top class in $\dl{R^{2^{n+1}}(n)}$ is $\dl{\zeta_{n+1}^{2^{n+1}},}$ which ``splits off" because $\dl{A(n)}$ acts trivially on it, yielding the desired class in $\dl{H^*A(n)}.$  The resulting short exact sequence of $A(n)$-modules yields a splitting of the Koszul-type resolution, so the corresponding class in Ext is a nonzero-divisor in the cohomology ring $\dl{H^*A(n).}$  See \cite{MS1} for details.

\smallskip
In \cite{Rog}, Rognes follows up on the work of Davis and Mahowald to fill in all the details for the computation of $\dl{\hbox{Ext}_{A(2)}\f2},$ the $\dl{E_2}$-term of the classical Adams {\SS} for $\dl{\pi_*(tmf).}$ One of his very useful observations (on page 44) is that the Davis-Mahowald {\SS} for the cocommutative Hopf algebras $A(n)$ is multiplicative, which we'll exploit in the proof of the conjecture.

\section{Strategy of the proof} \label{Strategy}

\bigskip
Now we outline the proof of the conjecture. First, as a notational shortcut, the phrase $\dl{y \in \Ext{A(k)}\f2}$ should be read as ``the class $\dl{y}$ is defined and nonzero in $\dl{\Ext{A(k)}\f2.}$"

We begin
with a simple observation:
\begin{lem} If $v_j^k \in \Ext{A(m)}\f2$, then $v_j^k \in \Ext{A(i)}\f2$
for $i = j, j+1, \dots ,m.$
\end{lem}

\noindent {\it Proof:} Simply notice that the restriction maps in
cohomology commute:

\begin{equation*}
\xymatrix{
   \Ext{A(m)}\f2 \ar[r]^{restr} \ar[d]^{restr} & \Ext{E(m)}\f2 \ar[d]^{restr} \\
     \Ext{A(i)}\f2 \ar[r]^{restr}& \Ext{E(i)}\f2. \\}
 \end{equation*}

We will prove the conjecture by induction, in a manner that might be paraphrased as follows: Assume that the conjecture  ``works for $\dl{v_{n-1},}$" then prove that it must also ``work for $\dl{v_n}$."

To begin, we assume that $\dl{v_{n-1}^{2^{n+1}}}$ is defined and nonzero in $\dl{\Ext{A(2n-1)}\f2.}$  As a (quite relevant) aside, we could actually assume as a ``base case" that
$\dl{v_{n-1}^{2^{n}}}$ is nonzero in $\dl{\Ext{A(2n-2)}\f2,}$ which is isomorphic to the $\dl{E_1^{0,*,*}}$-term of the Davis-Mahowald {\SS} for $\dl{\Ext{A(2n-1)}\f2.}$  Because $\dl{\Ext{A(2n-1)}\f2}$ is isomorphic to $\dl{\Ext{A}\f2}$ in this range, this class must be killed in the Davis-Mahowald {\SS,} and the only possible differential is
\[
d_1^{DM}(v_{n-1}^{2^n}) = v_{n-2}^{2^n} h_{2n-1}.
   \]
By the work of Rognes (\cite{Rog}, p 44), the cocommutativity of each $A(r)$ implies that the Davis-Mahowald spectral sequence is multiplicative, so
$\dl{d_1(v_{n-1}^{2^{n+1}})=0,}$ since $\dl{v_{n-1}^{2^{n}}}$ is a cycle.  Note that the bottom cell of $\dl{R^2(2n-1)}$ is in dimension $\dl{2^{2n+1},}$ too high to be the target of a $\dl{d_2^{DM}.}$  Since the higher $\dl{R^\sigma (2n-1)}$s are even more highly connected, we conclude that $\dl{v_{n-1}^{2^{n}}}$ is nonzero in $\dl{\Ext{A(2n-2)}\f2}$ actually implies that $\dl{v_{n-1}^{2^{n+1}}}$ is nonzero in $\dl{\Ext{A(2n-1)}\f2.}$  In any case, we will assume that $\dl{v_{n-1}^{2^{n+1}}}$ lives where we want it to.

Next, we use the class $\dl{v_{n-1}^{2^{n+1}} \in \Ext{A(2n-1)}\f2}$ to construct a nonzero product $\dl{v_{n-1}^{2^{n+1}} h_{2n+1} \in \Ext{A(2n)}(R^1(2n+1),\F_2),}$ the $\dl{E_1}$-term of the Davis-Mahowald {\SS} converging to $\dl{\Ext{A(2n+1)}\f2.}$ In this range \hfill \break $\dl{\Ext{A(2n+1)}\f2}$ is isomorphic to $\dl{\Ext{A}\f2,}$ so we can use a classical May {\SS} argument to show that $\dl{v_{n-1}^{2^{n+1}}h_{2n+1}}$ cannot live there.

We will conclude that there must be a Davis-Mahowald {\SS} class $\dl{ x \in E_1^{0,*,*} = \Ext{A(2n)}\f2}$ such that $\dl{ d_1^{DM}(x) = v_{n-1}^{2^{n+1}}h_{2n+1}.}$  We then use the multiplicativity of the {\SS} to rule out all other possibilities and conclude that this class $\dl{x}$ must in fact be the desired
$\dl{v_n^{2^{n+1}} \in \Ext{A(2n)}\f2.}$

To continue the proof from this base case, we will use the inductive hypothesis that $\dl{v_{n-1}^{2^{n+1+k}}}$ is nonzero in $\dl{\Ext{A(2n-1+k)}\f2 }$ to produce a nonzero class corresponding to $\dl{v_{n-1}^{2^{n+1+k}} h_{2n+1+k}}$ in the $\dl{E_1}$-term of the Davis-Mahowald spectral sequence for $\dl{\Ext{A(2n+1+k)}\f2 }$ that must be killed by a differential -- otherwise it would persist to $\dl{\Ext{A}\f2.}$  We then show that the only way that this class can be killed in the Davis-Mahowald {\SS} is if $\dl{v_n^{2^{n+1+k}}}$ is nonzero in $\dl{\Ext{A(2n+k)}\f2.}$

\section{An Example} \label{Example}

Given the complexity of the notation, it's best to work through an ``early" example before dealing with the details of the proof of the general case.  Here we detail the step from ``the conjecture is true for appropriate powers of $\dl{v_1}$" to ``it must be true for the appropriate powers of $\dl{v_2.}$"

Note that no power of $\dl{v_1}$ can be nonzero in $\dl{\Ext{A}\f2}$:  Such a class would be above the ``Adams edge," so there would be no possible targets for differentials on it in the classical Adams spectral sequence, yielding a nonzero nonnilpotent homotopy class in violation of Nishida's Theorem on the nilpotence of the stable homotopy ring of the sphere, per \cite{Adams} and \cite{Nishida}.

First, we observe that the conjecture holds for all appropriate powers of $\dl{v_1},$ by using well known calculations in the May {\SS.}  Note that $\dl{v_1^2}$ is represented in the May {\SS} by the class $\dl{b_{2,0},}
$ where we use the more modern notation from \cite{Ravbook}.

 This is as good a place as any to address concerns about whether we can identify these classes in a precise way in both the May and Davis-Mahowald spectral sequences. For arbitrarily chosen classes in $\dl{\Ext{A(k)}\f2,}$ there might be ambiguity in how one would choose representatives in these two spectral sequences.  For the particular classes we work with in this example (and in the proof of the general case), we can resolve this issue easily.  First, we note that $\dl{ h_k \in \Ext{A(k)}\f2}$ shows up in the Davis-Mahowald {\SS} in the $\sigma=1$ filtration, given by the class dual to $\dl{\zeta_1^{2^{k}} \in R^1(k),}$ which has an obvious counterpart in the May filtration of the cobar complex.  For the ``smaller" $\dl{h_j}$s in $\dl{\Ext{A(k)}\f2,}$ they show up in the Davis-Mahowald {\SS} in filtration 0 (alias $\dl{\Ext{A(k-1)}\f2)}$, so a simple induction ``back to" the $\dl{H^*(A(j))}$ case does the trick. The compatibility of the spectral sequence representatives for the $\dl{v_j^k}$s for $j<n$ is even easier to see: First, note that the ``last case" $\dl{v_n^{2^{n+1}} \in H^*A(n)}$ is easily resolved by the fact that there is only one May {\SS} generator in that bidegree, namely $\dl{b_{n-1,0}^{2^n}.}$  For the lower cases,
 the May filtration on $A(n)$ is compatible with the inclusion $\dl{E(n) \hookrightarrow A(n)},$ and the classes $\dl{v_j^k \in H^*A(n)}$ are actually defined in terms of the resulting restriction homomorphism.
\bigskip

Computations of Tangora, following May, in \cite{Tangora} show that
\[d_2^{May}(b_{2,0})  = b_{1,1} h_1 + b_{1,0} h_2 = h_1^3 + h_0^2 h_2.\]
   This differential is propagated by Nakamura's squaring operations \cite{Nak} (using the ``dual" versions $\dl{Sq_i}$ as in \cite{BEM}).  In particular, we will use the simplest case of these operations: $\dl{d_{2r}^{May}(Sq_0(x)) = Sq_1 (d_{r}^{May}(x)),}$ so the potential indeterminacy in the squaring operations is not an issue.

   We conclude that
  \begin{eqnarray*}
  d_4^{May}(b_{2,0}^2) & = & d_4^{May}(Sq_0 b_{2,0})\\
   & = & Sq_1(h_1^3 + h_0^2 h_2)\\
  & = & Sq_1 h_1^2 Sq_0 h_1 + Sq_1 h_1^2 Sq_0 h_1 + Sq_1h_0^2 Sq_0 h_2 + Sq_0 h_0^2 Sq_1 h_2\\
   & & \,\,\,\,\, \hbox{ (by the Cartan formula)}\\
  & = & h_1^4 h_1 + h_0^4 h_3 = h_0^4 h_3 \, \hbox{ (since $Sq_1 x^2 =0)$.}
  \end{eqnarray*}

We continue this process, to obtain
\[ d_{2^n}^{May}(b_{2,0}^{2^{n-1}}) = h_0^{2^n} h_{n+1}.\]
Tangora notes that in the range of his calculations, this May differential is necessary to truncate the $\dl{h_0}$-tower on $\dl{h_{n+1}}$ in $\dl{\Ext{A}\f2}$ at the desired height.

In fact, Tangora's reasoning can be used to see that $\dl{v_1^{2^{3+k}}}$ ``lives" as expected in $\dl{\Ext{A(1)}\f2}$ through $\dl{\Ext{A(3+k)}\f2,}$ meaning that the class is defined and nonzero in the cohomology of those $A(i)$s.  Precisely, if $\dl{v_1^{2^{3+k}}}$ is not a nonzero class in $\dl{\Ext{A(3+k)}\f2,}$ then there is no possible differential in the Davis-Mahowald {\SS} that can kill the class
$\dl{ h_0^{2^{3+k}} h_{4+k} \in E_1^{0,*,*}.}$ But if  $\dl{ h_0^{2^{3+k}} h_{4+k}}$ persists to a class in $\dl{\Ext{A(4+k)}\f2},$ it must also show up as a nonzero class in $\dl{\Ext{A}\f2,}$ by the Adams Approximation Theorem, contradicting the family of May differentials starting with Tangora's computation.  

An appropriate version of this reasoning is at the heart of the proof of the higher cases of the conjecture.

We will now show that the conjecture holding for all the appropriate powers of $\dl{v_1}$ implies that $\dl{v_2^{2^{3+k}}}$ must be nonzero in $\dl{\Ext{A(4+k)}\f2}$ for all $k \ge 0.$  The first nontrivial case is the following:

\centerline{$\dl{v_1^8 \in \Ext{A(3)}\f2}$ implies $\dl{v_2^8 \in \Ext{A(4)}\f2}$.}

 We first show that $\dl{v_1^8 \in \Ext{A(3)}\f2}$ implies $\dl{v_1^8h_5}$ is nonzero in the $\dl{E_1}$-term for $\dl{ \Ext{A(5)}\f2}.$   To see why, observe that the Davis-Mahowald {\SS} for $\dl{ \Ext{A(5)}\f2}$ has
\[
E_1^{\sigma,s,t} = \hbox{Ext}_{A(4)}^{s-\sigma,t}((R^\sigma(5))^*, \F_2),
\]
for all $\sigma \ge 0.$  We note that $\dl{E_1^{0,*,*}}$ is just $\dl{\hbox{Ext}_{A(4)}^{*,*}\f2},$  and the class \hfill \break
$\dl{h_5 \in \Ext{A(5)}\f2}$ is given by the bottom cell of $\dl{R^1(5)},$ via the map $\dl{\Sigma^{32} \F_2 \to R^1(5)}.$ We also note that $\dl{R^1(5)}$ is isomorphic to $\dl{\Sigma^{32} A(4)\sslash A(3)}$ as $\dl{A(4)_*}$-comodules, through dimension 62, which can see easily by comparing the cell diagrams.
More precisely, recall that for
$\dl{E_5 = E(\zeta_1^{32}, \zeta_2^{16},\dots,\zeta_{6}),}$
we have the following long exact sequence of $A(5)$-comodules, analogous to the Koszul resolution:
\[
0 \to \F_2 \to E_5 \otimes R^0(5) \to E_5 \otimes R^1(5) \to E_5 \otimes R^2(5) \cdots
\]
where $\dl{R^2(5)}$ begins in dimension 64.  Thus we have an $A(5)$-comodule isomorphism through degree 63 from $\dl{\F_2 \{\xi_1^{32}, \dots , \xi_6 \} \subset E_5 = E_5 \otimes R^0(5)}$ to \hfill \break
$\dl{\F_2 \{x_{32}, x_{48}, x_{56}, x_{60} ,x_{62},  x_{63} \} \subset E_5 \otimes R^1(5)}.$
So we know that \hfill  \break
$\dl{\Ext{A(4)}(R^1(5),\F_2) \cong \Ext{A(3)}\f2}$ in this range, and we must have a nonzero class $\dl{v_1^8h_5}$ corresponding to $\dl{`` \Sigma^{32} v_1^8" \in \Ext{A(4)}(R^1(5),\F_2)}. $

We can easily observe that our class $\dl{v_1^8h_5}$ cannot persist through the Davis-Mahowald {\SS} be a nonzero class in $\dl{\Ext{A(5)}\f2}$ (and hence in $\dl{\Ext{A}\f2}$) by looking at the ``next" family of May {\SS} differentials, starting with $\dl{b_{3,0}}$ corresponding to $\dl{v_2^2}$:
\begin{equation}\label{Maydiff}
d_2^{May}(b_{3,0})= b_{2,1}h_1 + b_{2,0}h_3.
\end{equation}
But $\dl{d_2^{May}b_{2,1} = h_2^3 + h_1^2 h_3,}$ so it follows that $\dl{b_{2,1} h_1 = h_1^3h_3.}$   Applying the Nakamura $\dl{Sq_0}$ twice to equation \ref{Maydiff}, we obtain
\[ d_8^{May}(b_{3,0}^4)=  b_{2,0}^4 h_5,\]
which shows that $\dl{v_1^8h_5}$ cannot be a nonzero class in $\dl{\Ext{A}\f2,}$ as we hoped.
\smallskip

A quick check of the dimensions of the modules $\dl{R^\sigma(5)}$ and the trajectories of the differentials in the Davis-Mahowald {\SS} shows that the only possible way to ``kill" the class $\dl{v_1^8h_5 \in \Ext{A(4)}(R^1(5),\F_2)}$ is to have a class \hfill  \break
 $\dl{x \in E_1^{0,8,56}= \Ext{A(4)}\f2}$ which could bear a $\dl{d_1.}$  We hope that the only candidate for such a class would be $\dl{v_2^8},$ IF we knew that
it persisted to $\dl{\Ext{A(4)}\f2}.$

 \smallskip

So we know that there must exist a class $\dl{ x \in \hbox{Ext}_{A(4)}^{8,56}\f2 = E_1^{0,8,56}}$ such that $\dl{d_1^{DM}(x) = v_1^8h_5.}$  We need to show that this class $x$ is exactly the desired $\dl{v_2^8 \in \Ext{A(4)}\f2.}$  First, we recall from the construction that the Davis-Mahowald {\SS} is multiplicative, since $\dl{A(k)}$ is a cocommutative Hopf algebra.  (See \cite{Rog} for details.)  Next, note that $\dl{v_1^8h_5}$ bears an $\dl{h_0}$-tower in $\dl{\Ext{A(4)}(R^1(5),\F_2),}$ so whatever classes and differentials are involved in killing it in the Davis-Mahowald {\SS} must also account for the entire $\dl{h_0}$-tower.  Recall that for a subalgebra $B$ of $A,$ the $\dl{h_0}$-towers in $\dl{\Ext{B}\f2 }$ are in one-to-one correspondence with the $\dl{Q_0}$-homology $\dl{H(A\sslash B;Q_0),}$ as detailed in \cite{Davis} \footnote{ Because this reference is not available online, we present an outline of Davis's proof. For any $A$-module $M,$ recall that the $\dl{Q_0}$-homology of $M$ is given by
$\dl{H(M;Q_0)= \dfrac{ker(Sq^1)}{Im(Sq^1)}.}$
Construct an epimorphism of $A$-modules $\dl{ M \mapright{\phi} N = (\bigoplus A) \oplus (\bigoplus A \sslash A(0)),}$
sending the $A$-module generators of $M$ to the generators of the first summand of $N$ and the generators of $\dl{H(M;Q_0)}$ to the generators of the second summand.  Davis observes that $L=ker(\phi)$ has zero $\dl{Q_0}$-homology, so that $\dl{\Ex{L,\F_2)}}$ is zero above the line $3s-t+6$ by the vanishing theorem of Adams \cite{Adams}.  Thus $\dl{\Ex{M,\F_2)} \cong \Ex{N,\F_2)}}$ above that line, as we hoped.
}.  In particular, since
\[
(A\sslash A(4))_* \cong \F_2[\xi_1^{32},\xi_2^{16}, \xi_3^8,\xi_4^4,\xi_5^2,\xi_6,\xi_7, \dots],
\]
we see that $\dl{H((A\sslash A(4))_*;Q_0)}$ is spanned as an $\dl{\F_2}$-vector space by $\dl{\chi \hbox{Sq}(32i,16j,8k,4l),}$ in dimensions $32i+48j+56k+60l,$ for $i,j,k,l \ge 0.$  So there's a unique tower in $\dl{\Ext{A(4)}\f2}$ in $t-s=48.$  Note that the restriction homomorphism \hfill \break
$\dl{restr: \Ext{A(4)}\f2 \to \Ext{E(4)}\f2}$ sends the elements in this tower to \hfill \break  $\dl{h_0^m v_2^8 \in \Ext{E(4)}\f2}$ for some $m$ (and all its $h_0$-multiples).  This does not show that $\dl{v_2^8}$ itself is present in $\dl{\Ext{A(4)}\f2,}$ but $\dl{h_0^k v_2^8}$ and all its $\dl{h_0}$-multiples must be nonzero for some (possibly large) $k.$

We might worry that the class we've detected, $\dl{x \in \Ext{A(4)}\f2},$ could be $\dl{h_0}$-torsion, but that the Davis-Mahowald {\SS} $\dl{d_1}$-differentials might still wipe out the whole tower on $\dl{v_1^8 h_5}$ in the manner given by the following diagram:

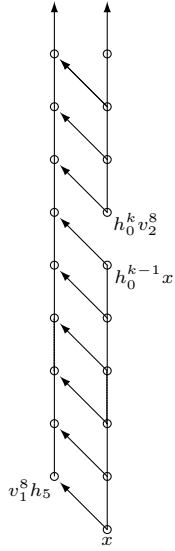
\begin{center}
\begin{picture}(230,230)
\put(111,39){\scsz \makebox(0,0)[t]{$v_1^8h_5$}}
\put(120,40){\elt}
\put(120,60){\elt}
\put(120,80){\elt}
\put(120,100){\elt}
\put(120,120){\elt}
\put(120,140){\elt}
\put(120,160){\elt}
\put(120,180){\elt}
\put(120,200){\elt}
\put(120,200){\tower}
\put(120,40){\hoa}
\put(120,60){\hoa}
\dottedline(120,80)(120,100)
\put(120,100){\hoa}
\put(120,120){\hoa}
\put(120,140){\hoa}
\put(120,160){\hoa}
\put(120,180){\hoa}
\put(140,20){\elt}
\put(140,17){\scsz \makebox(0,0)[t]{$x$}}
\put(140,20){\dia}
\put(140,40){\elt}
\put(140,40){\dia}
\put(140,60){\elt}
\put(140,60){\dia}
\put(140,80){\elt}
\put(140,80){\dia}
\put(140,100){\elt}
\put(140,100){\dia}
\put(140,120){\elt}
\put(154,120){\scsz \makebox(0,0)[t]{$h_0^{k-1}x$}}
\put(140,120){\dia}
\put(140,140){\elt}
\put(140,140){\dia}
\put(140,160){\elt}
\put(151,140){\scsz \makebox(0,0)[t]{$h_0^{k}v_2^8$}}
\put(140,160){\dia}
\put(140,180){\elt}
\put(140,180){\dia}
\put(140,200){\elt}
\put(140,180){\dia}
\put(140,200){\tower}
\put(140,20){\hoa}
\put(140,40){\hoa}
\dottedline(140,60)(140,80)
\put(140,80){\hoa}
\put(140,100){\hoa}
\put(140,140){\hoa}
\put(140,160){\hoa}
\put(140,180){\hoa}
\end{picture}
\end{center}

However, the multiplicativity in the Davis-Mahowald {\SS} prohibits this: If $x$ is $\dl{h_0}$-torsion and $\dl{d_1(x) = v_1^8 h_5,}$ then
\[
0 = d_1(h_0^L x) = h_0^L v_1^8h_5 \ne 0,
\]
for some large $L.$  We conclude that the class $x$ we detected in $\dl{\hbox{Ext}_{A(4)}^{8,56}\f2}$ must be exactly $\dl{v_2^8,}$ as we hoped.
\smallskip

This process can be continued to show that
$\dl{v_1^{16} \in \Ext{A(4)}\f2}$ implies \break $\dl{v_2^{16} \in \Ext{A(5)}\f2}$.
The argument is now easy to see: First, we observe that $\dl{R^1(6) \cong \Sigma^{64} A(5)\sslash A(4)}$ as $\dl{A(5)_*}$-comudules in the relevant range of dimensions, so
$\dl{v_1^{16}h_6 }$ must be a nonzero class in $\dl{\Ext{A(5)}(R^1(6),\F_2)},$ which contributes via the Davis-Mahowald {\SS} to $\dl{\Ext{A(6)}\f2}.$  If this class survives the Davis-Mahowald {\SS}, it must also live in $\dl{\Ext{A}\f2},$ by the Adams Approximation Theorem.  We use the May {\SS} differential \hfill  \break $\dl{d_{16}^{May}(b_{3,0}^8) = b_{2,0}^8 h_6}$ to see that  $\dl{v_1^{16}h_6}$ cannot persist to $\dl{\Ext{A}\f2.}$  The only way to kill $\dl{v_1^{16}h_6}$ in the Davis-Mahowald {\SS} for $\dl{\Ext{A(6)}\f2}$ is for some class $x$ to be nonzero in $\dl{E_1^{0,16,112} \cong \Ext{A(5)}\f2,}$ thus providing the source for the $\dl{d_1^{DMSS}}$ differential to kill $\dl{v_1^{16} h_6}.$  We see that $x$ must indeed be exactly $\dl{v_2^{16}}$ using the multiplicativity of the {\SS,} by examining the $\dl{Q_0}$-homology of $\dl{A\sslash A(5)}$ and seeing that there is only one $\dl{h_0}$-tower in $\dl{\hbox{Ext}_{A(5)}^{t-s=96}\f2,}$ corresponding to $\dl{h_0^k v_2^{16}.}$
\smallskip

The other powers of $\dl{v_2}$ are shown to live in the appropriate Exts similarly.

\section{Proof of the General Case } \label{Proof}

We prove that there is a nonzero class
\[ v_n^{2^{n+1+k}} \in \Ext{A(2n+k)}\f2,\]
by using induction on $n.$  We assume inductively that there is a nonzero element
\[ v_{n-1}^{2^{n+1+k}} \in \Ext{A(2n-1+k)}\f2.\]
We look closely at the Davis-Mahowald {\SS} for $\dl{\Ext{A(2n+1+k)}\f2,}$ focusing particularly on
$\dl{E_1^{1,*,*} \cong \Ext{A(2n+k)}(R^1(2n+1+k),\F_2).}$   The class \hfill \break
$\dl{h_{2n+1+k} \in \Ext{A(2n+1+k)}\f2}$ is detected in the $\dl{E_1^{1,*,*}}$-term by the bottom cell of $\dl{R^1(2n+1+k)}.$  More precisely, $\dl{h_{2n+1+k}}$ is given by applying $\dl{\Ext{A(2n+k)}(\,-\, ,\F_2)}$ to the homomorphism
$\dl{\Sigma^{2^{2n+1+k}} \F_2 \to R^1(2n+1+k).}$  We know that \hfill  \break
$\dl{ R^1(2n+1+k) \cong \Sigma^{2^{2n+1+k}} A(2n+k) \sslash  A(2n-1+k)}$ as $\dl{A(2n+k)_*}$-comodules through dimension $\dl{2^{2n+2+k}-2,}$ by examining the cell diagrams or by the following argument:  Observe, as in the example, that for
$\dl{E_r = E(\zeta_1^{2^r}, \zeta_2^{2^{r-1}},\dots,\zeta_{r+1}),}$
we have the following long exact sequence of $A(r)$-comodules, analogous to the Koszul resolution:
\[
0 \to \F_2 \to E_r \otimes R^0(r) \to E_r \otimes R^1(r) \to E_r \otimes R^2(r) \cdots
\]
where $\dl{R^2(r)}$ begins in dimension $\dl{2^{r+1}}$.  Thus we have an $A(r)$-comodule isomorphism through degree $\dl{2^{r+1}-1}$ from $\dl{\F_2 \{\xi_1^{2^r}, \dots , \xi_{2^{r+1}-1} \} \subset E_r = E_r \otimes R^0(r)}$ to \hfill \break
$\dl{\F_2 \{x_{2^r}, \dots,  x_{2^{r+1}-1} \} \subset E_r \otimes R^1(r)}.$

We conclude, then, that $\dl{v_{n-1}^{2^{n+1+k}} h_{2n+1+k}}$ is nonzero in the $\dl{E_1}$-term of the Davis-Mahowald {\SS} for $\dl{H^*A(2n+1+k).}$ If this class survived the {\SS} to the cohomology of
 {$A(2n+1+k),$} then the Adams Approximation Theorem tells us that it would be a nonzero class in $\dl{\Ext{A}\f2.}$  But in the May {\SS} for the cohomology of $A,$ we have the family of differentials starting with $\dl{d_2^{May}(b_{n+1,0}) = b_{n,1} h_{1} + b_{n,0} h_{n+1}.}$   We observe that \hfill  \break $\dl{d_2^{May}(b_{n,1}) = b_{n-1,2}h_2 + b_{n-1,1}h_{n+1}}$ and $\dl{d_2^{May}(b_{n,1}h_1) = 0 + b_{n-2,1}h_1h_{n+1},}$ so the ``extra" class in the $\dl{d_2}$ on $\dl{b_{n+1,0}}$ is dead by the $\dl{E_3}$-term.  We propagate this differential using Nakamura's squaring operations to obtain
\[ d_{2^{n+k}}^{May} (b_{n+1,0}^{2^n}) = v_{n-1}^{2^{n+1+k}} h_{2n+1+k}.
\]
So the class $\dl{v_{n-1}^{2^{n+1+k}} h_{2n+1+k}}$ cannot survive to $\dl{\Ext{A}\f2}$ (and hence, to \break  \hfill $\dl{\Ext{A(2n+k+1)}\f2,}$), and there must exist a class $\dl{ x \in \Ext{A(2n+k)}\f2 \cong E_1^{0,*,*}}$ such that
$\dl{d_1^{DM}(x) = v_{n-1}^{2^{n+1+k}} h_{2n+1+k}.}$

We show that the class $x$ must be exactly $\dl{v_n^{2^{n+1+k}}}$ by examining the $\dl{Q_0}$-homology of $\dl{(A\sslash A(2n+k))_* \cong \F_2[\xi_1^{2^{n+1+k}}, \xi_2^{2^{n+k}}, \dots ]}$ and observing the presence of a unique class in dimension $\dl{ 2^{n+k} \times 3},$ which must map under the restriction homomorphism to the class $\dl{v_n^{2^{n+1+k}} \in \Ext{E(2n+k)}\f2.}$  Thus $\dl{\Ext{A(2n+k)}\f2}$ contains an $\dl{h_0}$-tower starting at some $\dl{h_0}$ multiple of $\dl{v_n^{2^{n+1+k}}.}$  As in the $n=2$ example above, the $\dl{h_0}$-linearity of the Davis-Mahowald differentials demonstrates that the class $x$  must be exactly $\dl{v_n^{2^{n+1+k}},}$ as we wished.

\bibliographystyle{amsplain}

\end{document}